\documentclass[11pt,a4wide,oneside,reqno]{amsart}
\usepackage{amssymb,amsmath,url,a4wide,amsthm}
\usepackage[utf8]{inputenc}
\usepackage[left=3cm,right=3cm,top=3cm,bottom=2cm]{geometry}
\usepackage{enumitem}

\newcommand{\citep}[1]{\cite{#1}}

\newcommand{\F}{\mathbb{F}}
\newcommand{\N}{\mathbb{N}}

\newcommand{\R}{\mathbb{R}}

\newcommand{\e}{\varepsilon}
\newcommand{\inv}{^{-1}}
\newcommand{\lv}{\left\vert}
\newcommand{\rv}{\right\vert}
\newcommand{\lp}{\left(}
\newcommand{\rp}{\right)}
\newcommand{\calC}{\mathcal{C}}
\newcommand{\dlin}{d^\mathrm{lin}}

\renewcommand{\mod}{\;\mathrm{mod}\;}

\setlist[enumerate]{label=(\roman*)}

\usepackage[breaklinks=true]{hyperref}
\usepackage{cleveref}
\crefformat{footnote}{#2\footnotemark[#1]#3}

\newlist{algenum}{enumerate}{1}
\setlist[algenum]{label=Step \arabic*., ref=\arabic*}
\crefname{algenumi}{Substep}{substeps}
\Crefname{algenumi}{Step}{Steps}

\numberwithin{equation}{section}
\usepackage{thmtools}
\newtheorem{theorem}{Theorem}[section]

\newtheorem{lemma}[theorem]{Lemma}
\newtheorem{conjecture}[theorem]{Conjecture}
\newtheorem{proposition}[theorem]{Proposition}
\newtheorem{corollary}[theorem]{Corollary}
\newtheorem{question}[theorem]{Question}
\theoremstyle{definition}

\newtheorem{example}[theorem]{Example}

\declaretheoremstyle[%
  spaceabove=0pt,%
  spacebelow=0pt,%
  headfont=\normalfont\itshape,%
  postheadspace=1em,%
  qed=\qedsymbol%
]{mystyle}

\usepackage[url=false,doi=true,maxnames=4,giveninits=true,isbn=false]{biblatex}
\usepackage{xpatch}
\xpatchbibmacro{volume+number+eid}{%
  \setunit*{\adddot}%
}{%
}{}{}
\DeclareFieldFormat[article]{number}{\mkbibparens{#1}}
\renewbibmacro{in:}{%
  \ifentrytype{article}{}{\printtext{\bibstring{in}\intitlepunct}}}
\bibliography{biblio}

\title{Upper bounds for linear graph codes}

\author{Leo Versteegen}
\thanks{The author is grateful to be funded by Trinity College of the University of Cambridge through the Trinity External Researcher Studentship.}
\address{Department of Pure Mathematics and Mathematical Statistics, Centre for Mathematical Sciences, Wilberforce Road, Cambridge CB3 0WB, United Kingdom}
\email{lvv23@dpmms.cam.ac.uk}

\begin{document}
\begin{abstract}
A linear graph code is a family $\mathcal{C}$ of graphs on $n$ vertices with the property that the symmetric difference of the edge sets of any two graphs in $\mathcal{C}$ is also the edge set of a graph in $\mathcal{C}$. In this article, we investigate the maximal size of a linear graph code that does not contain a copy of a fixed graph $H$. 
In particular, we show that if $H$ has an even number of edges, the size of the code is $O(2^{\binom{n}{2}}/\log n)$, making progress on a question of Alon. Furthermore, we show that for almost all graphs $H$ with an even number of edges, there exists $\varepsilon_H>0$ such that the size of a linear graph code without a copy of $H$ is at most $2^{\binom{n}{2}}/n^{\varepsilon_H}$.
\end{abstract}
\maketitle
\section{Introduction}
The set of graphs with vertex set $[n]$ can be made into a vector space over $\F_2$ by defining the sum of two graphs $G_1=([n],E_1)$ and $G_2=([n],E_2)$ as the graph on $[n]$ whose edge set is the symmetric difference of $E_1$ and $E_2$. We denote this vector space by $\F_2^{K_n}$, and refer to subsets $\calC\subset \F_2^{K_n}$ as \emph{graph codes}. If $\calC$ is a subspace of $\F_2^{K_n}$, then we call it a \emph{linear graph code}. 
The focused study of graph codes was recently initiated by Alon et al. \cite{alon2023graphcodes,alon_et_al_structured} although graph codes had previously made implicit appearances elsewhere (see \cite{ellis_triangle_intersecting} for example).

In \cite{alon2023graphcodes}, Alon asked about the maximal size of a graph code such that the difference of no two graphs in the code is a copy of a fixed graph $H$. In this context, we say that $G$ is a \emph{copy} of $H$ if $[n]$ contains a subset $S$ such that $G[S]$ is isomorphic to $H$ and all vertices in $[n]\setminus S$ are isolated in $G$. For a fixed graph $H$, Alon defines
\begin{equation*}
d_H(n)=\frac{1}{2^{\binom{n}{2}}} \max\{\lv \calC\rv : \text{$\calC\subset \F_2^{K_n}$ and no sum of two elements in $\calC$ is a copy of $H$}\}.
\end{equation*}
Considering the code $\calC\subset \F_2^n$ that consists of all graphs on $[n]$ with an even number of edges, it is easy to see that $d_H(n)=1/2$ for all graphs $H$ with an odd number of edges. Alon also estimated the asymptotic behaviour of $d_H(n)$ for a few classes of graphs with strong symmetry properties such as matchings and stars, but beyond those classes, very little is known about the behaviour of $d_H(n)$. In particular, the following question of Alon is still open.

\begin{question}\label{question:K_4_general}
Is $d_{K_4}(n)=o(1)$?
\end{question}

Alon showed that the behaviour of $d_H(n)$ is closely linked to the following problem of a Ramsey-theoretic flavor. We say that an edge coloring $\chi$ of $K_n$ \emph{admits an even-chromatic copy} of $H$ if there exists an embedding $f\colon H \rightarrow K_n$ so that each color of $\chi$ is found at an even number of edges of $f(H)$. What is the minimum number $r_H(n)$ of colors of an edge coloring of $K_n$ that does not admit an even-chromatic copy of $H$? 

We first discuss how the quantity $r_H(n)$ relates to $d_H(n)$. Let $e(H)$ be the number of edges of $H$ and let $r=r_H(n)$. Suppose we are given a coloring $\chi\colon E(K_n)\rightarrow [r]$ that does not admit an even-chromatic copy of $H$. Suppose further that for some $t\in \N$, there exist $v_1,\ldots,v_r\in \F_2^t$ such that no sum of at most $e(H)$ of the $v_i$ is zero. Given such a family of vectors, the kernel of the unique linear map $\phi\colon \F_2^{K_n}\rightarrow \F_2^t$ that maps each edge $e\in E(K_n)$ to $v_{\chi(e)}$ does not contain a copy of $H$, which implies that $d_H(n)\geq 2^{-t}$. Alon remarked that one can obtain vectors $v_1,\ldots,v_r$ as above from the columns of a parity check matrix of a BCH code with designed distance $e(H)+1$ and length $t=\lceil (e(H)+1)/2\rceil \lceil \log_2 (r+1)\rceil$. Thus, if $r_H(n)$ is sufficiently large, then $d_H(n)\geq r_H(n)^{-(e(H)/2+1)}$.

The edge coloring of $K_n$ presented in \cite{cameron_K4_4_colors}, which is based on constructions from \cite{conlon_Kp_coloring,mubayi_Kp_small_coloring}, does not admit an even-chromatic copy of $K_4$ and uses (asymptotically) fewer than $\exp(\log^c n)$ colors for some $c<1$.\footnote{Throughout this article, $\log$ denotes the natural logarithm and $\log_b$ denotes the logarithm to base $b$.} After \cite{alon2023graphcodes} first appeared on the arXiv, Alon was informed of the existence of this coloring by Hunter and Mubayi and concluded in a revised version of his article that the following holds.

\begin{theorem}[Alon]\label{thm:K4-lower-bound}
There exists $c<1$ such that $r_{K_4}(n)\leq\exp(\log^c n)=n^{o(1)}$, and thus, $d_{K_4}(n)\geq n^{-o(1)}$.
\end{theorem}

The coloring in \cite{cameron_K4_4_colors} was not constructed with this purpose in mind, and it turns out that by simplifying it, we are able to prove \Cref{thm:K4-lower-bound} with $c=1/2$ (see \Cref{proposition:K_4-power2}).

To the best of the author's knowledge, all graph codes that have been used to establish lower bounds on $d_H(n)$ for any $H$ to date have been linear, so it is possible that linear graph codes are optimal in general. Regardless of whether this is true, Alon suggested that linear graph codes should be investigated in their own right, as is done for their counterparts in classical coding theory. Analogously to $d_H(n)$, we define

\begin{equation*}
\dlin_H(n)=\frac{1}{2^{\binom{n}{2}}}  \max\{\lv \calC\rv : \text{$\calC\subset \F_2^{K_n}$ is a subspace and does not contain a copy of $H$}\}.
\end{equation*}

Obviously, $d_H(n)\geq \dlin_H(n)$, and by the discussion above, $\dlin_H(n)\geq r_H(n)^{-(e(H)+1)/2}$. What is more, if $M$ is a parity-check matrix\footnote{We will repeatedly abuse language by pretending that for $N=\binom{n}{2}$, the spaces $\F_2^{K_n}$ and $\F_2^{N}$ are one and the same when formally, they are only isomorphic. Thus, a parity check matrix $M$ of $\calC\leq \F_2^{K_n}$ is a binary matrix with $\binom{n}{2}$ columns such that the kernel of $M$ is the isomorphic image of $\calC$ in $\F_2^{N}$.} of a linear code $\calC$ with co-dimension $t$ containing no copy of $H$, then the edge coloring that assigns to each edge $e$ the color $Me\in \F_2^{t}$ admits no even-chromatic copy of $H$. Thus,
\begin{align}\label{eq:coloring-linear}
\dlin_H(n)\leq \frac{1}{r_H(n)}.
\end{align}
As Alon observed, if $e(H)$ is even, then we can combine the bound $\dlin_H(n)\leq r_H(n)\inv$ with a multi-color Ramsey theorem to show that $\dlin_H(n)=O(\log\log n/\log n)$, which is the best known upper bound for general $H$. Alon suggested that it should be possible to improve this bound. Our first result confirms that one can at least remove the factor $\log\log n$.

\begin{theorem}\label{thm:general-bound}
There exists $c>0$ such that for all graphs $H$ with an even number of edges, $r_H(n)\geq c\log n$, and thus
\begin{align}\label{eq:general-bound}
\dlin_H(n)\leq (c\log n)^{-1}
\end{align}
for all sufficiently large $n$.
\end{theorem}

We suspect that \eqref{eq:general-bound} can be improved further, but due to the fact that $\dlin_{K_4}(n)\geq \exp(-\log^c n)$ with $c<1$, we cannot hope to improve it beyond a quasipolynomial bound in general. Remarkably, the case $H=K_4$ is a genuine bottleneck for \eqref{eq:general-bound} in the following sense. 

\begin{proposition}\label{thm:K4-bottleneck}
For all graphs $H$ with an even number of edges, 
\begin{align*}
r_H(n)\geq r_{K_4}\lp ne^{-O(\log^{3/4}n)} \rp.
\end{align*}
\end{proposition}

\Cref{thm:general-bound} is a corollary of \Cref{thm:K4-bottleneck} and a result of Fox and Sudakov \cite{fox_K4_lower_bound}, which asserts that \eqref{eq:general-bound} holds for $H=K_4$. Note that by \Cref{thm:K4-bottleneck}, any improvement to \eqref{eq:general-bound} for $H=K_4$ by at least a constant factor automatically yields an improvement on \eqref{eq:general-bound} in general. What is more, if it should turn out to be true that $r_{K_4}(n)\geq \exp(\log^{c}n)$ for some $c>0$, then we would obtain that $r_H(n)\geq \exp((1-o(1))\log^c n)$ for all $H$ with an even number of edges.

Setting aside $K_4$, it is not difficult to see that \eqref{eq:general-bound} can be improved for large classes of graphs. For instance, by the Kővári–Sós–Turán Theorem \cite{kovari_zarankiewicz}, for each bipartite graph $H$, there exists $c_H>0$ such that $r_H(n)\geq n^{c_H}$. This observation can be generalized to all graphs with an \emph{even decomposition}, which is a sequence of sets $V(H) =V_0\supset V_1 \supset \ldots \supset V_k=\emptyset$ such that for all $i\in [k]$, the set $V_{i-1}\setminus V_i$ is independent and the number of edges in $H$ between $V_{i-1}\setminus V_i$ and $V_i$ is even.

\begin{theorem}\label{prop:even-decomp-bound}
For every graph $H$ with an even composition, there exists $c_H>0$ such that $r_H(n)\geq n^{c_H}$, and thus $\dlin_H(n)\leq n^{-c_H}$. More precisely, if $H$ has $v(H)\geq 3$ vertices, then $r_H(n)\geq (1-o(1))n^{1/(v(H)-2)}$.
\end{theorem}

Clearly, $K_4$ does not have an even decomposition. On the other hand, it is also easy to see that every graph $H$ with an even number of edges that does not contain $K_4$ has an even decomposition. In fact, we will show that almost all graphs have an even decomposition.

\begin{theorem}\label{thm:even-decomp-prob}
The proportion of graphs on $n$ vertices with an even number of edges but without an even decomposition is at most $\exp(-\Omega(\sqrt{\log n}))$.
\end{theorem}

Together with \Cref{prop:even-decomp-bound}, this has the following immediate corollary.

\begin{corollary}
For asymptotically almost all graphs $H$ with an even number of edges, $r_H(n)\geq (1-o(1))n^{1/(v(H)-2)}$, and thus $\dlin_H(n)\leq (1+o(1))n^{-1/(v(H)-2)}$.
\end{corollary}

The remainder of the article is structured as follows. In \Cref{section:even-decomp}, we consider graphs with even decompositions, proving \Cref{prop:even-decomp-bound} and \Cref{thm:even-decomp-prob}. In \Cref{section:general-case}, we turn to the general case and prove \Cref{thm:general-bound} and \Cref{thm:K4-bottleneck}. Finally, in \Cref{section:concluding} we will make additional remarks on a few questions that arise from our results.

\section{Even decompositions}\label{section:even-decomp}
We first prove \Cref{prop:even-decomp-bound}. Its qualitative part follows immediately from repeated applications of the following lemma.

\begin{lemma}\label{lemma:decomp-step}
Suppose that a graph $H$ with an even number of edges contains an independent set $I$ of size $t$ so that for $F=H-I$, there exists $c_F>0$ such that $r_F(n)\geq (1-o(1))n^{c_F}$. Then for $c_H=c_F/(1+tc_F)$, we have $r_H(n)\geq (1-o(1))n^{c_H}$.
\end{lemma}
\begin{proof}
Let $\e>0$, $n\in \N$, and $r<(1-\e)n^{c_H}$, and let $\chi$ be an edge coloring of $K_n$ with $r$ colors. We consider the spanning subgraph $G\subset K_n$ whose edges have the most common color under $\chi$. By the pigeonhole principle, $G$ has at least $\binom{n}{2}/r$ edges, meaning that the average degree of $G$ is at least $n^{1-c_H}+t$ if $n$ is sufficiently large. We choose now a set $T\subset V(G)$ of size $t$ uniformly at random. The average size of the common neighbourhood $N(T)$ of $T$ is
\begin{align}\label{eq:dep-choice-avg-nbhd}
\frac{1}{\binom{n}{t}} \sum_{w\in G} \binom{d(w)}{t},
\end{align}
where we make the convention that $\binom{k}{t}=0$ for $k<t$. We extend the domain of the binomial coefficient to $\R$ by defining
\begin{align*}
\binom{\cdot}{t}\colon \R\rightarrow \R_{\geq 0} \qquad x\mapsto \binom{x}{t}=\begin{cases}
\frac{x(x-1)\cdots(x-t+1)}{t!} &\text{if $x\geq t-1$,}\\
0&\text{otherwise.}
\end{cases}
\end{align*} 
Note that this function is convex, and thus by Jensen's inequality, \eqref{eq:dep-choice-avg-nbhd} is at least
\begin{align*}
\frac{n}{\binom{n}{t}} \binom{\sum_{w\in G}d(w)/n}{t} \geq \frac{n}{\binom{n}{t}} \binom{n^{1-c_H}+t}{t}> n^{1-tc_H}.
\end{align*}
Thus, there exists a set $T\subset V(G)$ such that $\lv N(T)\rv> n^{1-tc_H}$. However, that means that
$\lv N(T)\rv^{c_F} >n^{c_H}>r$. Therefore if $n$ is sufficiently large, $\chi$ admits an even-chromatic copy of $H-I$ in $N(T)$, which we can extend to an even-chromatic copy of $H$ by mapping vertices from $I$ into $T$ using an arbitrary bijection. Since $\e$ may be chosen arbitrarily small, the Lemma follows.
\end{proof}

To complete the proof of \Cref{prop:even-decomp-bound}, it remains only to discuss the bound on the coefficient $c_H$.

\begin{proof}[Proof of \Cref{prop:even-decomp-bound}]
Suppose that a graph $H$ with $v(H)\geq 3$ has an even decomposition $V_0,\ldots, V_k$ and that $r_F(n)\geq (1-o(1))n^{1/(v(F)-2)}$ for all graphs $F$ with $3\leq v(F)\leq v(H)-1$. Let $I=V_0\setminus V_1$ and $t=\lv I\rv$. Without loss of generality, we may assume that $I$ has no non-empty subset $J$ such that $e(J,V_0\setminus J)$ is even, and in particular, that $t\leq 2$. Else, we may simply refine the decomposition until this is true. 
Consider now $F=H-I$. If $v(F)\in \{1,2\}$, then $H$ has at most four vertices, and $H$ must be the empty graph, a matching with two edges, or a path of length 2. In these cases, it is easy to see that $r_H(n)\geq (1-o(1))n$, leaving nothing to show.

We may therefore assume that $v(F)$ is at least 3, and because $V_1,\ldots,V_k$ is an even decomposition of $F$, we conclude that $r_F(n)\geq (1-o(1))n^{1/((v(F)-2)}$. By \Cref{lemma:decomp-step}, we have $r_H(n)\geq (1-o(1))n^{c_H}$ with
\begin{align*}
c_H=\frac{\frac{1}{v(F)-2}}{1+\frac{t}{v(F)-2}}=\frac{1}{v(F)-2+t}=\frac{1}{v(H)-2},
\end{align*}
as desired.
\end{proof}

We now turn our attention to \Cref{thm:even-decomp-prob}. When constructing an even decomposition $V_0,\ldots,V_k$ of a random graph $H$, we think of each set $V_i$ as being the set of \emph{remaining vertices} after the vertices in $V_0\setminus V_i$ have been \emph{removed}. We also refer to $H[V_i]$ as the \emph{remaining graph}. One simple strategy to construct an even decomposition is to greedily remove vertices with an even degree (in the remaining graph). As long as many vertices are left, it is unlikely that one cannot find such a vertex, but with positive probability, the last few vertices form a clique, at which point there is no way to continue. 

In the proof of \Cref{thm:even-decomp-prob}, we avoid this problem by setting aside an independent set $X\subset V(H)$ at the start of the algorithm. Then, we remove the lion's share of vertices greedily as suggested above, and for every remaining vertex $v\in V(H)\setminus X$ with odd degree, we can find with high probability another vertex $w$ with odd degree in $X\setminus N(v)$ so that we can remove $v$ and $w$ together. In the end, the only remaining vertices are in $X$ and thus isolated.

Of course, if we are not careful, then the probabilistic events which guide this process may not be independent. For example, conditioned on the event that $X$ is the largest independent set of $H$, the probability of every edge outside of $X$ to be present is slightly larger than $1/2$, skewing the probability of degrees being even. Stochastic dependence of events is the main difficulty in converting the sketch above into a proof, and it motivates the additional steps in the actual algorithm.

\begin{proof}[Proof of \Cref{thm:even-decomp-prob}]
For $n\in \N$, let $p=\lceil 2\log_{4/3} n \rceil$, $m=\lfloor n/(p+1)\rfloor$ and $q=\lfloor \sqrt{\log_2 m}/2\rfloor$. We partition $[n]$ into two sets $A$ and $B$ such that $\lv A\rv = m$ and $\lv B\rv$ is between $pm$ and $p(m+1)$. We assign to each element $v\in A$ a set $P_v$ such that all of these sets are disjoint. We also fix an arbitrary subset $Q_v\subset P_v$ of size $q$ for each $v\in A$. We now describe an algorithm that attempts to find an even decomposition of an arbitrary graph $H$. The success of some of the steps depends on $H$, and if one of the step fails, we abort and consider the entire algorithm failed. Later, we will show that if $H$ is drawn uniformly at random from all graphs on $[n]$ (which is equivalent to sampling each edge independently with probability $1/2$), then the algorithm succeeds with high probability. 
\begin{algenum}[leftmargin=2.0cm]
\item\label{step:idp-set} Fix an independent set  $X\subset A$ of size $q^2$.
\item\label{step:fix-a} Fix an arbitrary vertex $a\in A\setminus X$.
\item\label{step:clean-A} Let $k=\lv A\setminus (X\cup \{a\})\rv$ and execute for each $i\in [k]$ \cref{step:clean-A}.$i$ as follows. Let $v$ be the smallest\footnote{\label{footnote:smallest}With respect to the natural ordering of $[n]\subset \N$.} remaining vertex in $A\setminus (X\cup \{a\})$. Query the parities of the degrees of $v$ and all $w\in P_v$ in the remaining graph. Query also which edges between $v$ and $P_v$ are present. If the degree of $v$ is even, remove $v$.  Otherwise, if there is a vertex $w\in P_v$ that is not adjacent to $v$ and has odd degree, remove $v$ and $w$. If no such vertex $w$ exists, then abort the algorithm.
\item\label{step:clean-a} Proceed as in \Cref{step:clean-A} for $a$ but use $Q_a$ instead of $P_a$.
\item\label{step:main} Denote by $C$ all remaining vertices in $B\setminus Q_a$ at the start of this step. We execute a sequence of substeps \ref{step:main}.$i$ for $i\in \{1,\ldots,\lv C\rv-q\}$, removing one vertex from $C$ in each substep until only $q$ of them are left. In \cref{step:main}.$i$, we query the degrees of all remaining vertices from $C$ in the remaining graph. If none of them have even degree, then abort the algorithm. Otherwise, remove the smallest\cref{footnote:smallest} vertex with even degree.
\item\label{step:partition-X} Let $Y$ be the set of all remaining vertices from $B$, including those in $Q_a$. Partition $X$ into sets of size $q$ using their order in $\N$ and assign to each element of $v\in Y$ one of the partition sets $R_v$.
\item\label{step:fix-b} Fix an arbitrary element $b\in Y$.
\item\label{step:clean-B} 
\noindent Let $k=\lv Y\setminus \{b\}\rv$ and execute for each $i\in [k]$ \cref{step:clean-B}.$i$ as follows. Let $v$ be the smallest\cref{footnote:smallest} remaining vertex in $Y\setminus \{b\}$. Query the parities of the degrees of $v$ and all $w\in R_v$ in the remaining graph. Query also which edges between $v$ and $R_v$ are present. If the degree of $v$ is even, remove $v$.  Otherwise, if there is a vertex $w\in R_v$ that is not adjacent to $v$ and has odd degree, remove $v$ and $w$. If no such vertex $w$ exists, then abort the algorithm.
\item\label{step:final} The only remaining vertices are in $X\cup \{b\}$. If $H$ has an even number of edges, then the degree of $b$ in the remaining graph is even. In this case, we remove $b$ and then the remaining independent set, which shows that $H$ has an even decomposition.
\end{algenum} 

Each step of the algorithm consists of a combination of the following five actions:

\begin{itemize}
\item Labelling a vertex (or a set of vertices) for a later purpose.
\item Checking whether a certain edge is present.
\item Checking whether a certain vertex has even degree in a certain induced subgraph of $H$.
\item Removing a vertex.
\item Aborting the algorithm.
\end{itemize}

The second and third type of action can be understood as taking the inner product $\langle F,H\rangle$ of some \emph{query graph} $F\in \F_2^{K_n}$ with $H$, where the inner product of two graphs in $\F_2^{K_n}$ is the parity of the size of their intersection. These inner products determine the flow of the algorithm completely. In particular, if two graphs give the same result for all queries up to some step $\sigma$, the algorithm removes the same vertices from both of them and the same edges will be queried for both of them in the next step. In this situation, we say that the two graphs are \emph{equivalent up to $\sigma$}.

Let $S\subset \F_2^{K_n}$ be the set of graphs that are subject to a given (sub-)step $\sigma$, let $C$ be an equivalence class of equivalence up to $\sigma$ on $S$, and let $F_1,\ldots,F_k$ be the graphs that are used to query sets of edges up to $\sigma$ for the graphs in $C$. Let further $G_1,\ldots,G_{l}\in \F_2^{K_n}$ and $Z\subset \F_2^{l}$ be such that $\sigma$ fails for a graph $H\in C$ if and only if the vector $(\langle G_i,H\rangle)_{i\in [l]}$ is in $Z$. If we assume for the moment that the graphs $F_1,\ldots,F_k,G_1,\ldots,G_{l}$ are linearly independent, then $\sigma$ will fail for $\lv C\rv\cdot \lv Z \rv/2^{l}$ graphs in $C$. Looking at the steps of the algorithm, it is clear that the quantity $\lv Z \rv/2^l$ does not depend on the class $C$. Denoting it by $\varepsilon_\sigma$, we conclude that $\sigma$ fails for at most $\varepsilon_\sigma\lv S\rv\leq \varepsilon_\sigma 2^{\binom{n}{2}}$ graphs.

Two things remain to be shown. First, that for each graph $H$, the graphs $F_1,\ldots,F_k$ that are used to query parities throughout the algorithm are indeed linearly independent. Secondly, that the sum of $\varepsilon_\sigma$ over all (sub-)steps $\sigma$ is of order $\exp(-\Omega(\sqrt{\log n}))$.

Suppose now that for some graph $H$ the algorithm uses query graphs that are not linearly independent, meaning that for some query graphs $F_1,\ldots,F_k$ we have
\begin{align}\label{eq:zero-sum}
G:=F_1+\ldots+F_k=0.
\end{align}
In particular, for every edge $e\in \F_2^{K_n}$, we must have $\langle e,G\rangle=\langle e,F_1\rangle+\ldots+\langle e,F_k\rangle=0$. We now go through the different steps and argue that if one of the vectors $F_1,\ldots,F_k$ had been used to query the parity of a quantity during that step, then $F_1+\ldots+F_k$ could not be 0, which leads to a contradiction.

First, no $F_i$ can query the degree parity of any $v\in A\setminus X$ (including $a$) during \Cref{step:clean-A} or \Cref{step:clean-a}. This is because for any $u\in P_a\setminus Q_a$, we would have $\langle uv,G\rangle=\langle uv,F_i\rangle=1$, as the edge $uv$ is not part of any other query. Using this, we may conclude that no $F_i$ queries the presence of an edge within $A$ as in \Cref{step:idp-set} either. Next, for any $v\in A\setminus (X\cup\{a\})$, no $F_i$ queries the degree parity of any $w\in P_v$ during \Cref{step:clean-A} as otherwise we have $\langle wa,G\rangle=1$. From this, we conclude in turn that no $F_i$ queries the edges between $v$ and $P_v$ for any $v\in A\setminus (X\cup\{a\})$ as in \Cref{step:clean-A}.

Suppose now there is $j$ such that $F_j$ queries the degree parity of a vertex during \Cref{step:main}. We may assume without loss of generality that no graph among $F_1,\ldots,F_k$ queries the parity of a degree during an earlier substep of \Cref{step:main}. Let $v$ be the vertex such that $F_j$ queries the parity of the degree of $v$ during \cref{step:main}.$i$, and note that the algorithm cannot have failed during \cref{step:main}.$i$. Indeed. if this were the case, we would have $\langle vw,G\rangle=\langle vw,F_j\rangle=1$ for all $w\in X$, , which contradicts \eqref{eq:zero-sum}. Letting $u$ be the vertex that is removed during \cref{step:main}.$i$, we may assume that $u=v$, as otherwise $\langle uv,G\rangle=1$. However, this implies once more that $\langle vw,G\rangle=\langle vw,F_j\rangle=1$ for all $w\in X$.

Let now $w\in Q_a$ and let $u$ be the first vertex that was removed during \Cref{step:main}. Since none of $F_1,\ldots,F_k$ queries the degree parity of $u$, none of them can query the degree parity of $w$ in \Cref{step:clean-a} either because that would imply $\langle uw,G\rangle=1$. We may now also conclude that none of the graphs query the presence of an edge between $a$ and $Q_a$ during \Cref{step:clean-a}.

Analogously to what we said about \Cref{step:clean-A} and \Cref{step:clean-a}, we can argue that no $F_j$ queries the degree parity of $v\in Y$ or of $w\in R_v$ for $v\in Y\setminus \{b\}$ during \Cref{step:clean-B} or \Cref{step:final}.

We now have to show that the sum of the quantities $\varepsilon_\sigma$ over all (sub-)steps $\sigma$ is small. We think of $\varepsilon_\sigma$ as the probability that the edges under investigation during $\sigma$ do not have suitable parities if each edge is sampled independently with probability $1/2$. Regarding \Cref{step:idp-set}, it is well known that a random graph on $m$ vertices has an independent set of size larger than $\log_2 m$ with high probability. We gain little by bounding this probability too carefully, so for simplicity, we bound the probability that $A$ has no independent set of size $2q^2$ by $e^{-m^{3/2}}<e^{-n}$ for large $m$ (see e.g. \cite{frieze_random_graphs} for a bound which implies this one).

The failure probability for each substep of \Cref{step:clean-A} is $(3/4)^p/2$ because we only fail to remove $v$ if the degree of $v$ is odd and for each $w\in P_v$, either the degree of $w$ is even or $vw$ is present.
Similarly, each substep of \Cref{step:clean-B} fails with probability $(3/4)^q/2$ and \Cref{step:clean-a} fails with probability $(3/4)^q/2$. For each $i\in [\lv C\rv-q]$, \Cref{step:main}.$i$ fails with probability $2^{i-1-\lv C\rv}$, meaning that we can bound the probability that \Cref{step:main} fails overall by $2^{-q}$. \Cref{step:final} fails if the degree of $b$ is odd, which happens with probability $1/2$. Thus the proportion of graphs with an even number of edges but without an even decomposition is less than
\begin{align*}
e^{-n}+\frac{3^pn}{4^p}+\frac{3^q}{4^q}+2^{-q}+\frac{3^q\cdot 2q}{4^q}.
\end{align*}
Since $p>2\log_{4/3} n$, we can bound the above by $(2+o(1))3^qq^2/4^q$. By definition, $q=\lfloor \sqrt{\log_2 m}/2\rfloor$, which is $(1-o(1))\sqrt{\log_2 n}$. Thus, 
\begin{align*}
(2+o(1))3^qq^2/4^q<2\exp(-(\log(4/3)-o(1))\sqrt{\log_2 n}+\log\log_2 n),
\end{align*}
which is of order $\exp(-\Omega(\sqrt{\log n}))$, as desired.
\end{proof}

\section{The general case}\label{section:general-case}
The main ingredient in the proof of \Cref{thm:K4-bottleneck} will be \Cref{lemma:superpolynomial-chrom-bound}, which we will use to decompose graphs that contain $K_4$.

\begin{lemma}\label{lemma:superpolynomial-chrom-bound}
Let $H$ be a graph whose vertex set can be partitioned into two non-empty sets $A_1$ and $A_2$ such that $e(A_1,A_2)$ is even. Let further $k=\lv A_1\rv$,  $r_1=r_{H[A_1]}$, and $r_2=r_{H[A_2]}$. Then for all sufficiently large $n$,
\begin{align}\label{eq:maxmin}
r_H(n)\geq \max_{m\in [n]} \min \left\{r_1(m),r_2(m),\exp\lp \frac{(\log n-\log m)^2}{3k\log n}\rp\right\}
\end{align}
\end{lemma}

In order to prove \Cref{lemma:superpolynomial-chrom-bound}, we will make use of the well-known combinatorial technique of dependent random choice in the form of the following lemma, which was first presented by Fox and Sudakov \cite{fox_dependent_random_choice}. Recall that for a graph $G$ and a set $S\subset V(G)$, we denote by $N_G(S)$ the common neighbourhood of $S$, i.e., the vertices in $G$ that are adjacent to every $v\in S$. 

\begin{lemma}[Dependent random choice]\label{lemma:dep-choice}
Let $G$ be a graph on $n$ vertices with at least $\alpha n^2/2$ edges for some $\alpha>0$. Let further $u,k,m,t\in \N$ such that
\begin{align}\label{eq:dep-choice}
n\alpha ^{t}-\binom{n}{k}\left({\frac {m}{n}}\right)^{t}\geq u.
\end{align}
There exists a set $U\subset V(G)$ such that $\lv U \rv\geq u$ and for all $S\subset U$ with $\lv S\rv=k$, $\lv N_G(S)\rv\geq m$.
\end{lemma}

\begin{proof}[Proof of \Cref{lemma:superpolynomial-chrom-bound}]
For $n\in \N$, let $r=r(n)$ be the maximum on the right hand side of \eqref{eq:maxmin}, and let $m=m(n)$ be an integer that achieves this maximum. To prove the claim, it is enough to show that for all sufficiently large $n$, any coloring $\chi$ of $K_n$ with less than $r$ colors admits an even-chromatic copy of $H$. Given such a coloring, we define $G$ as the spanning subgraph of $K_n$, whose edges are those that have the most frequent color under $\chi$. In particular, $G$ has more than $n^2/2r$ edges. We now want to apply \Cref{lemma:dep-choice} with $k$, $m$, $u=m$, $\alpha=1/r$, and 
\begin{align*}
t=\left\lceil \frac{k\log n}{\log n - \log m - \log r} \right\rceil.
\end{align*}
By \eqref{eq:maxmin}, we have $\log r \leq (\log n-\log m)/3k$, whence $t$ is well-defined, positive, and 
\begin{align*}
t\leq (1+o(1))\frac{k\log n}{(1-1/3k)(\log n -\log m)}.
\end{align*}
Hence,
\begin{align*}
n\alpha^t&\geq \exp\lp \log n - (1+o(1))\frac{k\log n \log r}{(1-1/3k)(\log n -\log m)}\rp \\
&\geq \exp \lp \log n - (1+o(1))\frac{\log n - \log m}{3-1/k}\rp,
\end{align*}
Since the claim is trivially true unless $\log n -\log m(n) \rightarrow \infty$, we may assume that the above is at least $2m$ if $n$ is sufficiently large. At the same time, we have that 
\begin{align*}
\frac{\binom{n}{k}\lp \frac{m}{n}\rp^t}{n\alpha^t}&<\exp\lp (k-1)\log n + t(\log m - \log n + \log r) \rp\\
&<\exp\lp (k-1)\log n - (k-o(1))\log n\rp.
\end{align*}
Since the above is $o(1)$, \eqref{eq:dep-choice} is satisfied for all sufficiently large $n$, and we can find a set $U\subset [n]$ larger than $m$ such that for all subsets $S\subset U$ of size $\lv A_1\rv$, there is a set $M\subset [n]\setminus U$ of size $m$ such that all edges between $S$ and $M$ have the same color.

Since $r\geq r_1(m)$, there is a set $S_1\subset U$ such that $\chi$ admits an even-chromatic copy of $H[A_1]$ on $S_1$. Furthermore, because $r\geq r_2(m)$ as well, the $N_G(S_1)$ contains a set $S_2$ such that $\chi$ admits an even-chromatic copy of $H[A_2]$ on $S_2$. Since all edges between $S_1$ and $S_2$ are in $G$ and $e(A_1,A_2)$ is even, we know that $\chi$ admits an even-chromatic copy of $H$ on $S_1\cup S_2$.
\end{proof}

To prove \Cref{thm:K4-bottleneck}, we require the following lower bound for $r_{K_4}(n)$. The coloring in the proof of \Cref{proposition:K_4-power2} is a simplified version of the coloring constructed by Cameron and Heath \cite{cameron_K4_4_colors}.

\begin{proposition}\label{proposition:K_4-power2}
$r_{K_4}(n)\leq \exp(O(\sqrt{\log n}))$
\end{proposition}
\begin{proof}
Let $d=\lceil \sqrt{\log_2 n}\rceil$ and $m=2^d$, and note that $m^d\geq n$. Denoting the lexicographical order on $V=[m]^d$ by $<$, we define $\psi\colon V^{(2)}\rightarrow \{-1,0,1\}^d$, which maps $vw$ to $\psi(vw)$, where for all $i\in [d]$,
\begin{align*}
\psi_i(vw)=\begin{cases}
1&\text{if  $(v< w$ and $v_i < w_i)$ or $(v>w$ and $v_i>w_i)$,}\\
-1&\text{if  $(v< w$ and $v_i > w_i)$ or $(v>w$ and $v_i<w_i)$,}\\
0&\text{if  $v_i=w_i$.}
\end{cases}
\end{align*}
Furthermore, for each pair $vw \in V^{(2)}$, we let $\delta(vw)$ be the index of the first coordinate in which $v$ and $w$ differ. With these definitions in hand, we define an edge coloring $\chi$ of the complete graph on the vertex set $V$ by setting $\chi(vw)=(\psi(vw),\{v_{\delta(vw)},w_{\delta(vw)}\})$. Note that $\chi$ uses less than $3^dm^2=\exp(O(\sqrt{\log n}))$ colors, and since we can embed $K_n$ into the complete graph on $V$, it remains only to show that $\chi$ admits no even-chromatic copy of $K_4$.

Suppose towards a contradiction that there are vertices $v,w,x,y\in V$ which form an even-chromatic copy of $K_4$. Let $i$ indicate the first coordinate in which any of the four vertices differ, and note that $\{v_i,w_i,x_i,y_i\}$ contains only two values. Indeed, if $v_i$, $w_i$, $x_i$ and $y_i$ were all different, all six edges would have different colors, and if, for example, $v_i=w_i$ but $v_i$, $x_i$, $y_i$ were distinct, the edges $v_iw_i$, $v_ix_i$, $v_iy_i$, and $x_iy_i$ would all have different colors. Neither can be the case.

We may therefore assume without loss of generality that $v_i=w_i$ and $x_i=y_i$. Clearly, the only edge that can have the same color as $v_iw_i$ is $x_iy_i$, as for the other four edges we have $\psi_i(e)=1$.
Thus, if we let $j$ indicate the first coordinate on which $v$ and $w$ differ, then $j$ must also indicate the first coordinate on which $x$ and $y$ differ, and the set $\{v_j,w_j,x_j,y_j\}$ must have size two. Again without loss of generality, $v_i<x_i$, $v_j<w_j$, and $v_j=x_j$. Consequently, $w_j=y_j$, and $v<w<x<y$. However, this means that the color $\chi(vy)$ occurs only once among the six edges since $\psi_i(vy)=\psi_j(vy)=1$ but $\psi_i(vw)=\psi_i(xy)=0$, $\psi_j(vx)=\psi_j(wy)=0$ and $\psi_j(wx)=-1$. We arrive at a contradiction, and conclude that $\chi$ admits no even-chromatic copy of $K_4$.
\end{proof}

We are now ready to prove \Cref{thm:K4-bottleneck}.

\begin{proof}[Proof of \Cref{thm:K4-bottleneck}]
Let $H$ be a graph with an even number of edges such that the claim holds for all graphs with fewer vertices. If $\lv V(H)\rv\leq 1$, then $r_H(n)=\infty$, whence we may assume that $H$ has at least two vertices. By distinguishing whether $H$ is complete or not, it is easy to see that $V(H)$ can be partitioned into two non-empty sets $A_1$ and $A_2$ such that the numbers of edges in $A_1$, in $A_2$, and between $A_1$ and $A_2$ are all even. Letting $k=\lv A_1\rv$, \Cref{lemma:superpolynomial-chrom-bound} tells us that
\begin{align}\label{eq:minmax-specific}
r_H(n)\geq \max_{m\in [n]} \min \left\{r_{H[A_1]}(m),r_{H[A_2]}(m),\exp\lp \frac{(\log n-\log m)^2}{3k\log n}\rp\right\}.
\end{align}
Omitting rounding signs for the sake of clarity, we let $m= n\exp(-C\log^{3/4} n)$ for some $C>0$. To complete the proof, it is enough to show that if $C$ is sufficiently large, then all three terms on the right hand side of \eqref{eq:minmax-specific} are bounded from below by $r_{K_4}(n\exp(-C'\log^{3/4} n))$ for some $C'>0$. Since $H[A_1]$ and $H[A_2]$ have fewer vertices than $H$, there are $C_1,C_2>0$ so that for $i\in \{1,2\}$
\begin{align*}
r_{H[A_i]}(m)\geq r_{K_4}( m\exp(-C_i\log^{3/4} m))\geq r_{K_4}(n\exp(-(C+C_i)\log^{3/4} n)),
\end{align*}
which settles the first two terms in \eqref{eq:minmax-specific}. However, the third term is equal to $\exp(C^2\sqrt{\log n}/3k)$, and if $C$ is sufficiently large, then this is larger than $r_{K_4}(n)$ by \Cref{proposition:K_4-power2}, completing the proof.
\end{proof}

Now that we have proved \Cref{thm:K4-bottleneck}, \Cref{thm:general-bound} is an immediate corollary of the following lower bound on $r_{K_4}(n)$ due to Fox and Sudakov \cite{fox_K4_lower_bound}.

\begin{theorem}\label{thm:fox-K4-lower-bound}
For $k>2^{100}$, and $n\geq 2^{2000k}$, any edge coloring of $K_n$ with $k$ colors admits a copy of $K_4$ that is either monochromatic or the union of a monochromatic $C_4$ and a monochromatic matching.
\end{theorem}

\begin{proof}[Proof of \Cref{thm:general-bound}]
\Cref{thm:fox-K4-lower-bound} implies that $r_{K_4}(n)\geq c\log n$ for some $c>0$ and all sufficiently large $n$. But by \Cref{thm:K4-bottleneck}, if $H$ is a graph with an even number of edges, then
\begin{align*}
r_H(n)\geq r_{K_4}\lp ne^{-O(\log^{3/4}n)} \rp \geq c (\log n - O(\log^{3/4}n)).
\end{align*}
For sufficiently large $n$, this is at least $\frac{c}{2}\log n$.
\end{proof}

\section{Concluding remarks}\label{section:concluding}
Although we made progress in a few directions, many questions about the size of linear codes remain open. \Cref{thm:K4-bottleneck} shows that in a qualitative sense, $\dlin_{K_4}(n)$ is a lower bound for $\dlin_H(n)$ for all $H$ with an even number of edges. That being said, we are no closer to determining $\dlin_{K_4}(n)$. We suspect that the quasipolynomial lower bound for $\dlin_{K_4}(n)$ is closer to the truth than the logarithm upper bound. 

We have shown that graphs with even decompositions form a large class of graphs for which $\dlin_H(n)$ decreases at least as fast as $n^{-1/(v(H)-2)}$. Note that unless $H$ is empty, $\dlin_H(n)$ cannot possibly decrease faster than polynomially. Indeed, by sampling a parity-check matrix at random and employing the first moment method, we see that $\dlin_H(n)\geq n^{-v(H)}$. Using Lovász's local lemma, we can improve this to $\dlin_H(n)=\Omega(n^{-v(H)+2})$. Nevertheless, the proof of \Cref{prop:even-decomp-bound} is rather wasteful, and the bounds on the exponent $c_H$ can certainly be improved for many graphs.

We also believe that there are more graphs with even decompositions than \Cref{thm:even-decomp-prob} suggests. In fact, considering the extreme wastefulness of its proof, we conjecture that the following much stronger bound holds.

\begin{conjecture}
There exists $c>0$ such that the proportion of graphs on $n$ vertices with an even number of edges but without an even decomposition is at most $\exp(-cn^2)$.
\end{conjecture}

However, it cannot be case that the number of graphs without even decompositions is $2^{o(n^2)}$. Indeed, it is easy to see that if a graph is a union of a complete graph of size $2k+4$ and any other graph on $k$ vertices, then it has no even decomposition. Generally speaking, it seems to be easier to find even decompositions the sparser a graph is. At the same time, even very dense graphs can have an even decomposition. For example, it is easy to see that whenever $n$ is even, removing a matching from $K_n$ yields an evenly decomposable graph. We therefore propose a weighted version of \Cref{thm:even-decomp-prob}, which would bound the probability $\alpha_p(n)$ that the random graph $G(n,p)$ for $p=p(n)\neq 1/2$ lacks an even decomposition given that $G(n,p)$ has an even number of edges. The main obstacle in adapting the proof of \Cref{thm:even-decomp-prob} to the regime $p<1/2$ is that the relevant stochastic events cease to be independent. While this complicates things considerably, we are convinced that $\alpha_p(n)=o(1)$ if $p(n)\leq 1/2$. On the other hand, if $p(n)$ converges to 1 quickly enough, $\alpha_p(n)$ must converge to 1.

\begin{question}
For $n\in \N$, what is the minimal $p(n)$ such that $\alpha_p(n)\geq 1/2$?
\end{question}

Concerning graphs without even decompositions, Ge, Xu, and Zhang \cite{ge_xu_zhang_K5} and independently, Bennett, Heath, and Zerbib \cite{bennett_K5} recently showed that $r_{K_5}(n)$ also grows at most quasipolynomially. In \cite{ge_xu_zhang_K5}, the authors conjecture that $r_K(n)=n^{o(1)}$ holds for all cliques $K$. We think it is plausible that an even stronger statement is true.

\begin{conjecture}
For all graphs $H$ without an even decomposition, $r_H(n)=n^{o(1)}$.
\end{conjecture}

Finally, there is the question of how $\dlin_H(n)$ compares to $d_H(n)$. We are not aware of any graph for which these quantities are not equal, or even of a graph $H$ for which the largest known code $\calC$ such that $\calC- \calC$ contains no copy of $H$ is not linear. This begs the following question.

\begin{question}\label{question:d-dlin-equal}
Are $\dlin_H(n)$ and $d_H(n)$ (asymptotically) equal?
\end{question}

The only evidence we have against an affirmative answer to \Cref{question:d-dlin-equal} is the following example of a non-linear code $\calC$ that contains no triangle and has almost the optimal density of $1/2$.

\begin{example}
Let $\phi\colon \F_2^{K_n}\rightarrow \F_2^{K_n}$ be the unique linear map that maps each edge to its complement in $K_n$. If $n$ is 0 or 1 $\mod 4$, then $\phi$ is invertible (in fact, an involution), and letting $\mathcal{A}\subset\F_2^{K_n}$ be the set of graphs with at most $\binom{n}{2}/2-2$ edges, we can define $\calC=\phi\inv(\mathcal{A})$. This code has $(1-O(1/n))2^{\binom{n}{2}-1}$ elements, and we claim that $\calC-\calC$ contains no triangle. Indeed, suppose that for $G,H\in \calC$, the symmetric difference $G+H$ is a triangle. By definition of $\mathcal{A}$, we know that $\phi(G)$ and $\phi(H)$ have at most $\binom{n}{2}/2-2$ edges, meaning that $\phi(G+H)=\phi(G)+\phi(H)$ hast at most $\binom{n}{2}-4$ edges. But $\phi$ maps every graph with an odd number of edges to its complement, and in particular, it maps a triangle to a graph with $\binom{n}{2}-3$ edges, which yields a contradiction.
\end{example}

\printbibliography

\end{document}